\documentclass[journal]{IEEEtran}
\usepackage{graphicx} 
\usepackage{xcolor}
\usepackage{amssymb}
\usepackage{hyperref}

\usepackage{imports}

\usepackage{enumitem, tikz}

\usetikzlibrary{graphs}

\def\checkmark{\tikz\fill[scale=0.4](0,.35) -- (.25,0) -- (1,.7) -- (.25,.15) -- cycle;}

\begin{document}
\title{Locational Marginal Prices Obey DC Circuit Laws}

\author{Kyri Baker,~\IEEEmembership{Senior Member,~IEEE,} and Harsha Gangammanavar}


\maketitle

\begin{abstract}
Electricity markets often utilize the DC approximation of the AC power flow equations to facilitate solving an otherwise complex nonconvex optimization problem. These DC power flow equations have analogies to DC circuit laws such as Kirchhoff's Laws, resulting in an intuitive understanding of power flows under this model. Variables derived from the Lagrangian dual of the DC optimal power flow problem, such as locational marginal prices (LMPs) and congestion cost, are less intuitive without an understanding of optimization theory. Even with this understanding, LMP behavior, such as the conditions in which negative prices occur or the impact of individual congested lines on network-wide LMPs, remain somewhat mysterious. In this paper, we show that prices also obey DC circuit laws, which can help facilitate an intuitive understanding of their behavior and relationships throughout a network without explicitly understanding duality. In particular, prices can be modeled as voltages, and their differences can be modeled as flows, allowing for a physical interpretation of prices. This analogy also lends itself to the use of well-understood DC circuit concepts such as superposition and Kirchhoff's Laws, which can further facilitate a clearer understanding of price behavior.
\end{abstract}


\IEEEpeerreviewmaketitle

\section{Introduction}
Many independent system operators (ISOs) leverage the DC approximation of the AC power flow equations to clear power markets and determine locational marginal prices (LMPs) \cite{HistoryofOPF}. LMPs are typically said to represent ``the marginal cost of supplying an additional unit of power at a particular bus'' or some variation thereof. These values, in reality, do not actually represent the true cost to deliver power to a particular location, due to the physical approximations made within the DC OPF \cite{bakerOPF21}. As one example, the basic form of DC OPF does not include losses, or only includes loss estimations which are added to the LMP \cite{LineLossDC}. To be mathematically precise, LMPs are derived from the Lagrange multipliers (dual variables) corresponding to the power balance equation at each bus. 

Even with these simplifications, LMPs are generally challenging for many to understand from a technical perspective, even for some who work with them on a daily basis. While primal variables (voltages, power flows) can be intuitively understood and visualized (e.g. power flowing into a bus must equal power flowing out, current flows from a high to a lower voltage, etc.), prices exhibit somewhat strange behavior to many. For example, LMPs at adjacent buses can be wildly different, congestion in one portion of the grid can impact prices at buses that should seemingly be unaffected, negative prices sometimes occur for unintuitive reasons, and more. 

Thus, this paper aims to provide a more intuitive understanding of LMPs using concepts taught in high-school physics and introductory circuits courses. This allows for a more straightforward understanding of concepts such as how a single congested line impacts prices across the network (using superposition), how LMPs at adjacent buses are related, when negative prices occur, etc. Using this analogy, LMPs can be thought of as voltages and LMP differences as proportional to current flows between buses - facilitating the visualization of price behavior as physical quantities.

\begin{figure}[t!]
    \centering
    \includegraphics[width=1\columnwidth]{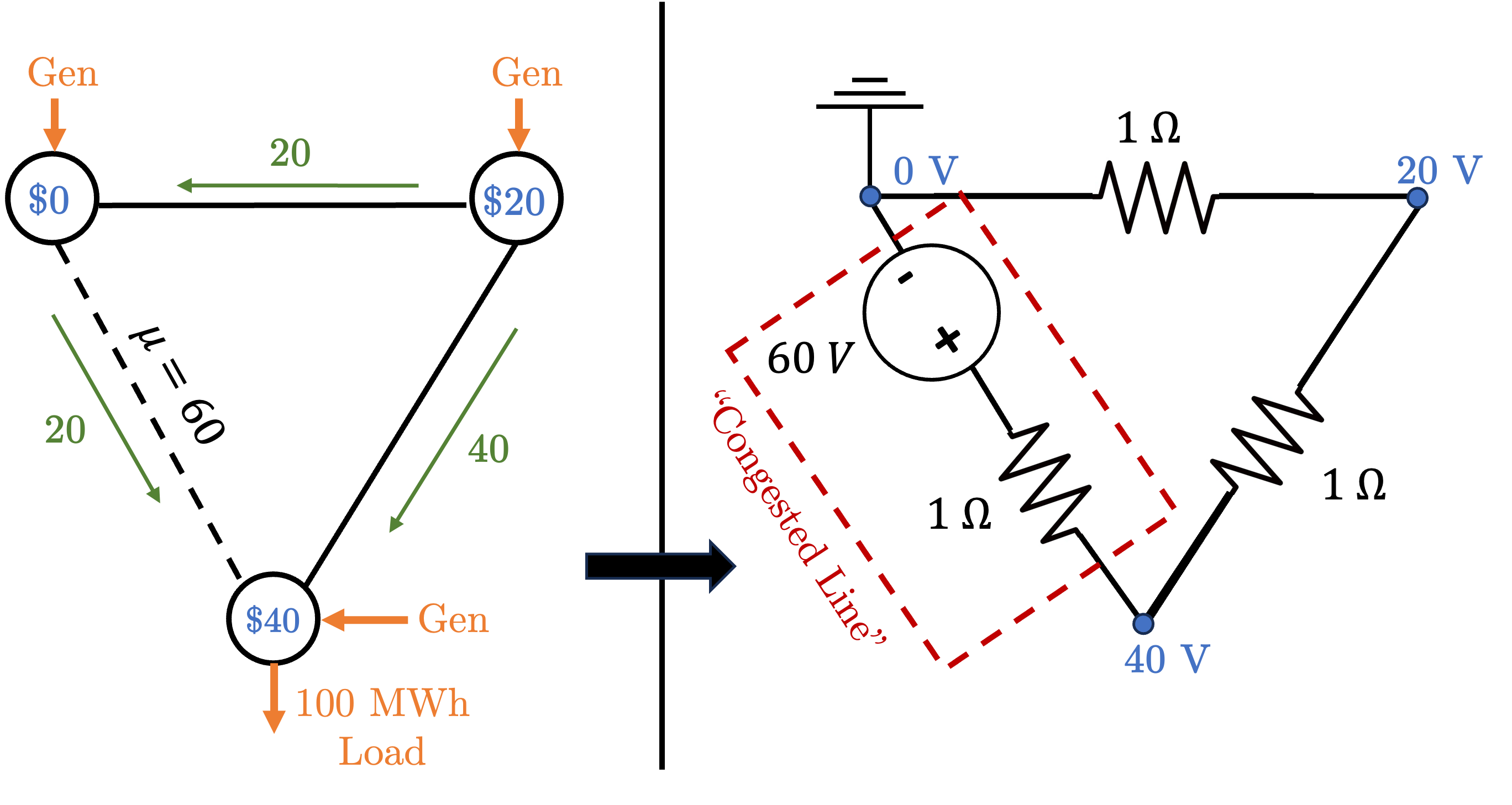}
    \caption{Left: A 3-bus network with 3 marginal generators and one load. LMPs are shown in blue at each bus, MW flows are shown between buses in green, and congestion (and its associated cost $\mu$) is shown with a dashed line. Right: The DC circuit equivalent: line susceptances become resistors, congestion becomes a voltage (or current) source, and ground is set at the node of the lowest cost marginal generator. The voltages in the circuit are analogous to LMPs in the DC OPF problem.}
    \label{fig:3bus_analogy}
\end{figure}

We will show the process of converting a power network's dual solution into an ``equivalent circuit'' in which DC circuit analysis can be performed. One example of this is shown in the left subfigure of Fig. \ref{fig:3bus_analogy}, for a three-bus network with three generators (of marginal costs \$0, \$20, and \$40) and one load at the bottom bus. All line susceptances are the same and equal to 1 pu, and the line between the two left-most buses is congested with a congestion cost of \$60. In the right subfigure, the analogous DC circuit is formed. Congested lines become independent voltage (in series with a resistor) or current (in parallel with a resistor) sources, lines become resistors (equal to the reciprocal of the original line susceptance), and LMPs are nodal voltages. The placement of ground, which impacts nodal voltages, is set to the node at which the cheapest marginal generator is located. If this generator does not have a cost of \$0, all voltages must be added to whatever marginal cost this is in order to recover the exact LMPs. More detail on the full analogy will be given in a later section.

DC circuit concepts have been used in education and research to solve the primal problem; for example, superposition can be used to expedite OPF solution times \cite{Li24}. Outside of the well-known DC circuit interpretation of the primal DC OPF problem, circuit analysis techniques have also been used to solve AC power flow problems as well more efficiently \cite{Jereminov19}, again focusing on creating equivalent circuit models for primal quantities. These techniques have also been used to solve challenging three-phase distribution power flow problems in a more robust manner \cite{Pandey19}. 

In this paper, we do not claim that the use of DC circuit solvers is necessarily faster than existing techniques to clear markets and determine LMPs, but rather that this interpretation of the dual problem lends itself to a more intuitive understanding of price behavior. Further, the use of DC circuit techniques such as superposition, Kirchhoff's Laws, Thevenin/Norton's theorems, etc. can allow for better conceptualization of the propagation of congestion on prices throughout the network, the distribution of price differences and how it relates to line susceptances and the difference in marginal generator costs, and the impact of individual congestions on LMPs network-wide, to name a few applications.

\section{Prices in the Dual Problem}
We first consider the DC OPF primal problem for a network of $n$ buses, which solves for the active power generation $\bp_g \in \reals^{n \times 1}$ and complex voltage angles at each bus $\btheta \in \reals^{n \times 1}$:
\begin{subequations} \label{eqn:DCOPF}
\begin{align} 
\min_{\bp_g, \btheta} ~~~&\bc^T \bp_{g} \\
\mathrm{s.t:~}& \bp_{g} + \bB\btheta = \bp_{d}, \label{eqn:DCPowerBalance}\\
& \underline{\bp}_d \leq \bp_d \leq \overline{\bp}_d,\\
& \underline{\bp}_g \leq \bp_g \leq \overline{\bp}_g.
\end{align}
\end{subequations}

\noindent for a network with $n$ buses and $l$ lines, we have the $n \times n$, symmetric ``DC bus admittance matrix'' $\bB$. Each off-diagonal entry $(i,j)$ in $\bB$ is $-b_{ij}$, or $-\frac{1}{x_{ij}}$. Each on-diagonal entry $(i,i)$ in $\bB$ is the sum of all of the $b_{ij}$'s corresponding to buses connected to bus $i$. For example, consider the three-bus example in Fig. \ref{fig:3bus_analogy}. The $\bB$ matrix would be:

\[
  \bB =
  \left[ {\begin{array}{ccccc}
    b_{12}+b_{13} & -b_{12} & -b_{13} \\
    -b_{12} & b_{12}+b_{23} & -b_{23} \\
    -b_{13} & -b_{23} & b_{13}+b_{23} \\
  \end{array} } \right]
\]

\noindent Note that the $\bB$ matrix is singular - this can be addressed by choosing a reference bus $i$ and angle $\theta_i = 0$. The reduced $\bB$ matrix, which removes the row and column corresponding to $\theta_i$ in this case, is typically not singular. This matrix is used, for example, in the ``shift factor'' formulation of the DC power flow equations, where the shift factors depend on the choice of reference bus. 

\subsection{The primal problem}

We can rewrite the DC OPF as a generic linear program: 
\begin{align} \label{eqn:DCOPF}
\min_{\bp, \btheta} ~~~&\bc^T \bp \tag{\bP} \\
\mathrm{s.t.~}& \bA \bp + \bB \btheta = \ba, \label{eq:KCL}\\
& \bC \bp \geq \bb,  \\
& \bD \btheta \geq \bd, 
\end{align}

\noindent where the equation satisfying (but not exactly equivalent to) Kirchhoff's law is \eqref{eq:KCL} - flows into a bus (e.g. from incoming line flows or generators) must be equal to flows out of a bus (e.g. to loads or outgoing line flows) \cite{Frank16}. In Kirchhoff's current law (KCL), these flows represent currents - in DC OPF, complex voltages, and currents are not modeled directly, but the flow balance concept is often used to explain active power balance as well.

For completeness, we can also consider loads that make bids into the market as well as the generators that make offers. Denote by $\bp \in \reals^{2n \times 1}$ a vector of power injections (though generation) or consumption (through controllable loads) at each bus, that is, $\bp = (\bp_g, -\bp_d)^\top$; and $\btheta \in \reals^{n \times 1}$ is a vector of voltage angles at each bus in the network.  $\ba \in \reals^{n \times 1}$ captures the uncontrollable loads at each bus (e.g., inflexible demands); $\bc \in \reals^{2n \times 1}$ represents the cost of a generator (to produce energy) or controllable load (to reduce energy) at each bus in the network. 

Matrix $\bA$ can be considered a $n \times 2n$ ``location matrix'' which contains a $1$ to denote generators and $-1$ to denote loads. We assume that a maximum of 1 generator and/or controllable load can be at a single bus (otherwise, the bus can be split into two, with a very small reactance between the two). E.g. for the 3-bus network in \ref{fig:3bus_analogy} with generators at each bus and a controllable load at bus 3, we have $\bp = [p_{g1}, p_{g2}, p_{g3}, 0, 0, p_{l3}]^T$, a (primal) objective of $c_1 p_{g1} + c_2 p_{g2} + c_3 p_{g3} - c_6 p_{l3}$, $\bc = [c_1, c_2, c_3, 0, 0, -c_6]^T$, $\ba = [0, 0, \bar{p}_{l3}]^T$, and $\bA = [\bI_n -\bI_n]$, where $\bI_n$ is the $n \times n$ identity matrix. matrices $\bD$ ($\bd$) and $\bC$ ($\bb$) represent constraint matrices (right-hand sides) relevant to the angle and power injection variables, respectively. Using the same example, assuming lower ($\underline{\bp}$, a vector of all lower bounds) and upper ($\overline{\bp}$, a vector of all upper bounds) bounds exist on all values in $\bp$, we have $\bC$ as a $4n \times 2n$ matrix and $\bb$ as a $4n \times 1$ vector. 

\subsection{The dual problem}

\noindent The linear programming dual of $\eqref{eqn:DCOPF}$ is given by:
\begin{align} \label{eqn:DCOPF_dual}
\max_{\substack{\blambda, \bmu, \bgamma}} ~~~&\ba^T \blambda + \bb^T \bmu + \bd^T \bgamma \tag{\bD} \\
\mathrm{s.t:~}& \bA^T \blambda + \bC^T \bgamma = \bc, \notag \\
& \bB^T \blambda + \bD^T \bmu = 0, \notag\\ 
& \bgamma, \bmu \geq 0. \notag
\end{align}
Where $\blambda$ are the dual variables associated with the equality constraints (power flow equations), $\bmu$ are the multipliers associated with the inequality constraints on $\btheta$ (line flow limits), and $\bgamma$ are the multipliers associated with the inequality constraints on $\bp$ (upper and lower generation limits). The optimality conditions of \eqref{eqn:DCOPF}-\eqref{eqn:DCOPF_dual} are thus:
\begin{subequations} \label{eqn:optimalityConditions}
\begin{align}
    & \bA \bp + \bB \btheta = \ba, \label{eqn:primalFeasibility}\\
    & \bA^T \blambda + \bC^T \bgamma = \bc, \label{eqn:dualFeasibility_part}\\
    & \bB^T \blambda + \bD^T \bmu = 0, \label{eqn:dualFeasibility_theta}\\
    & 0 \leq \bgamma \perp \bC \bp - \bb \geq 0, \label{eqn:complementary_part}\\
    & 0 \leq \bmu \perp \bD \btheta - \bd \geq 0, \label{eqn:complementary_theta}\\
    & \bgamma, \bmu \geq 0. \label{eqn:nonNeg_duals}
\end{align}
\end{subequations}

The dual variables (Lagrange multipliers) corresponding to the power balance equations (equality constraints) are the LMPs at each bus, $\blambda$, and variables $\bmu$ are representative of the price of congestion on each line.

\section{Conversion to a DC Circuit}
Interestingly, we show that in meshed networks, the behavior of the LMPs $\blambda$ and congestion prices $\bmu$ across the network can be modeled in a way akin to a linear DC circuit. As in the primal DC OPF problem \cite{Kocuk16}, these dual quantities obey Kirchhoff’s Voltage and Current Laws and other theorems such as superposition.

\subsection{Circuit Component Analogies}\label{sec:analogy}

To facilitate looking at LMP relationships through this lens, we first redefine dual concepts with their analogous linear circuit counterparts. For clarity, we use the term ``node'' to refer to nodes within a circuit, and the term ``bus'' to refer to buses within a power network. Lastly, a ``marginal generator'' is a generator $i$ with $\underline{p}_i < p^*_i < \overline{p}_i$, where $^*$ denotes the optimal value for that variable.

\begin{itemize}
    \item \textbf{Ground.} Ground is placed at the cheapest marginal generator, or the circuit node with the lowest nodal voltage. Note that although the concept of ground in a circuit is used as a reference point, this is distinctly different than the concept of a ``reference'' or ``slack bus'' in a power system (usually where $\theta_i$ o bus $i$ is fixed to zero).
    \item \textbf{Voltage.} The locational marginal price $\lambda_i$ at a bus $i$ is akin to the ``voltage'' at that bus. If the cheapest marginal generator does not have a marginal cost of $\$0$, the voltages at each node must be added to the cost of the cheapest marginal generator to recover the exact LMP magnitudes. 
    \item \textbf{Resistance.} The inverse of the susceptance $b_{ij}$ between buses $i$ and $j$, $\frac{1}{b_{ij}}$, is the ``resistance'' of the line between bus $i$ and $j$.
    \item \textbf{Current.} Thus, from Ohm's law, the ``current'' flowing between bus $i$ and $j$ is $b_{ij}(\lambda_i - \lambda_j)$.
    \item \textbf{Current Source.} Congested lines, that is, lines $ij$ with maximum allowable flow between bus $i$ and $j$, have $\mu_{ij} > 0$. This $\mu_{ij}$ acts as an ideal ``current source,'' parallel to the congested line, injecting a fixed amount of current from bus $i$ to $j$.
    \item \textbf{Voltage Source.} Alternatively to the current source analogy, a source transformation can be used to instead represent a current source in parallel with a resistor to a voltage source in series with a resistor.
\end{itemize}

Note that this analogy only holds if congestion is present in the network; otherwise no ``source'' will exist in the circuit. Without congestion, LMPs throughout a network should all be equal, and thus an analogy isn't necessary.  With these concepts, we can first start analyzing the applicability of Kirchhoff's Laws. 

\subsection{KCL and KVL}

\textbf{Kirchhoff's Current Law (KCL).} KCL, or conservation of charge, states that the total ``current'' flowing into a bus must equal the ``current'' flowing out of a bus. In the primal problem, KCL is represented by constraint $\bA\bp + \bB\btheta = \ba$. In the dual, using the definitions for current source, voltage, current, and resistance above, this follows from \eqref{eqn:dualFeasibility_theta}, where $\bB^T \blambda + \bD^T \bmu = 0$. Note the entries of $\bD$ are 0, -1, or 1.

\textbf{Kirchhoff's Voltage Law (KVL).} KVL, or conservation of energy, states that the sum of ``voltage drops'' around a loop must equal zero. In the primal problem, KVL holds around a loop by summing the difference in voltage angles $\btheta$ around a loop. In the dual, this is akin to summing the LMP differences around a loop, which would obviously sum to zero because each $\lambda_i$ would appear exactly once in the sum with a + sign and exactly once with a - sign.

\section{Case Study}

To illustrate these concepts, we use a meshed 7-bus example with congestion (otherwise all LMPs would be the same). This is a modified version of the 6-bus single cycle network from \cite{Pritchard10}, with an additional bus and two lines added, and an additional congested line.

In this example, all of the line susceptances are assumed to be equal to 1. Two lines are congested (indicated by dashed lines): the line between buses 1 and 6, and the line between buses 2 and 4. Blue numbers inside each circular bus symbol denote the LMP at that bus. Green ``current flows'' in the network are, again, the LMP differences between those buses multiplied by the line susceptance between those buses. Orange arrows into, or out of, a bus, denote generation (into) or demand (out of) that bus. In this example, there are 5 generators, two of which are producing at their upper limit, and three of which are marginal (at buses 1, 2, and 6). There is one load located at bus 4. 

\begin{figure}[t!]
    \centering
    \includegraphics[width=1\columnwidth]{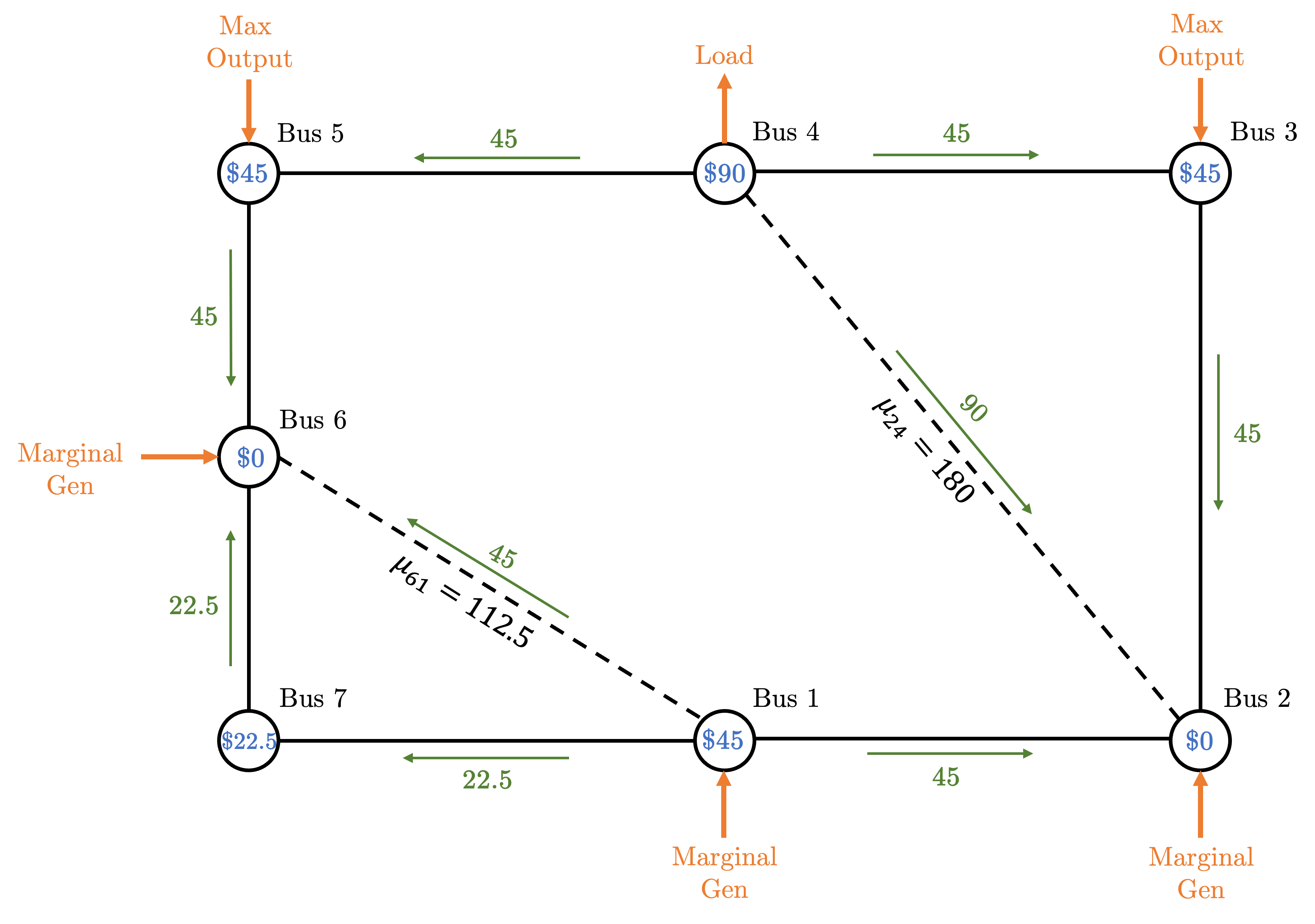}
    \caption{7-bus network with mesh topology. Dashed lines indicate congestion. Blue text indicates LMPs, and green text indicates LMP differences between buses multiplied by the susceptance of the line between those buses. Three generators are marginal (not at a lower or upper generation limit).}
    \label{fig:7bus}
\end{figure}

Additionally indicated on the figure is the congestion cost of the two lines that are congested. Recall that even though the ``flow'' of LMP is, say, from bus 1 to bus 6, the actual flow of power is from the cheaper generator at bus 6 towards bus 1. This is the reason that this multiplier is written $\mu_{61}$ rather than $\mu_{16}$. 

\subsection{Illustration of Kirchhoff's Laws}

As an illustrative example, we can perform a subset of KVLs around \emph{some} loops, recalling that LMPs are analogous to voltages.
\begin{align*}
    &\textit{Lower Left Loop: } (45-22.5) + (22.5-0) + (0-45) = 0 ~\checkmark\\
    &\textit{Center Loop: } (45-0) + (0-45) \\
    &\hspace{2cm}+ (45-90) + (90-0) + (0-45) = 0 ~\checkmark\\
    &\textit{Upper Right Loop: } (0-90) + (90-45) + (45-0) = 0 ~\checkmark\\
    &\textit{Outer Loop: } (0-45) + (45-22.5) + (22.5-0) \\&\hspace{0.7cm}+ (0-45) + (45-90) + (90-45) + (45-0) = 0 ~\checkmark
\end{align*}

KCLs can be written at each bus by summing the ``current'' into (left hand side) and exiting (right hand side) a bus, recalling that congested lines add an additional, parallel, current source equal to the value of the Lagrange multiplier ($\mu$) corresponding to the binding flow limits of that line.

\begin{align*}
    \textit{KCL at Node 1: } &112.5 = 45 + 22.5 + 45 ~\checkmark\\
    \textit{KCL at Node 2:  } &45 + 90 + 45 = 180 ~\checkmark\\
    \textit{KCL at Node 3:  } &45 = 45 ~\checkmark\\
    \textit{KCL at Node 4:  } &180 = 45 + 90 + 45~\checkmark\\
    \textit{KCL at Node 5:  } &45 = 45~\checkmark\\
    \textit{KCL at Node 6:  } &45 + 45 + 22.5 = 112.5~\checkmark\\
    \textit{KCL at Node 7:  } &22.5 = 22.5~\checkmark
\end{align*}

In this version of Kirchhoff's laws, primal variables like power injections and voltage angles do not appear at all.

\subsection{Price Differences as ``Flows''}
An improved understanding of how LMPs propagate from node to node can now be interpreted using the KCL concept. Consider the top-left bus in Fig. \ref{fig:7bus}. \$45 ``flows'' into that bus, and thus \$45 must ``flow out.'' This creates a \$45 difference between buses 4, 5, and 6. If the susceptance of line 4-5 or line 5-6 were different, this would result in a different distribution of LMPs, since the analogy of current is actually the LMP difference between two buses multiplied by the susceptance of the line between those buses. 

Just as a voltage or current source in one area of a circuit can still impact nodal voltages multiple ``hops'' away, this confirms that even though congestion does not exist in that area of the grid, congestion in other parts of the grid will still affect flows (e.g. LMP differences) in nodes in other areas of the circuit/grid. Similar intuition could be used to deduce how some grid-enhancing technologies such as flexible alternating current transmission system (or FACTS) devices can be used to adjust line parameters such as susceptance \cite{ADETOKUN2021e06461} to change LMP flows.

\subsection{Conversion to a circuit}

\begin{figure}[t!]
    \centering
    \includegraphics[width=0.8\columnwidth]{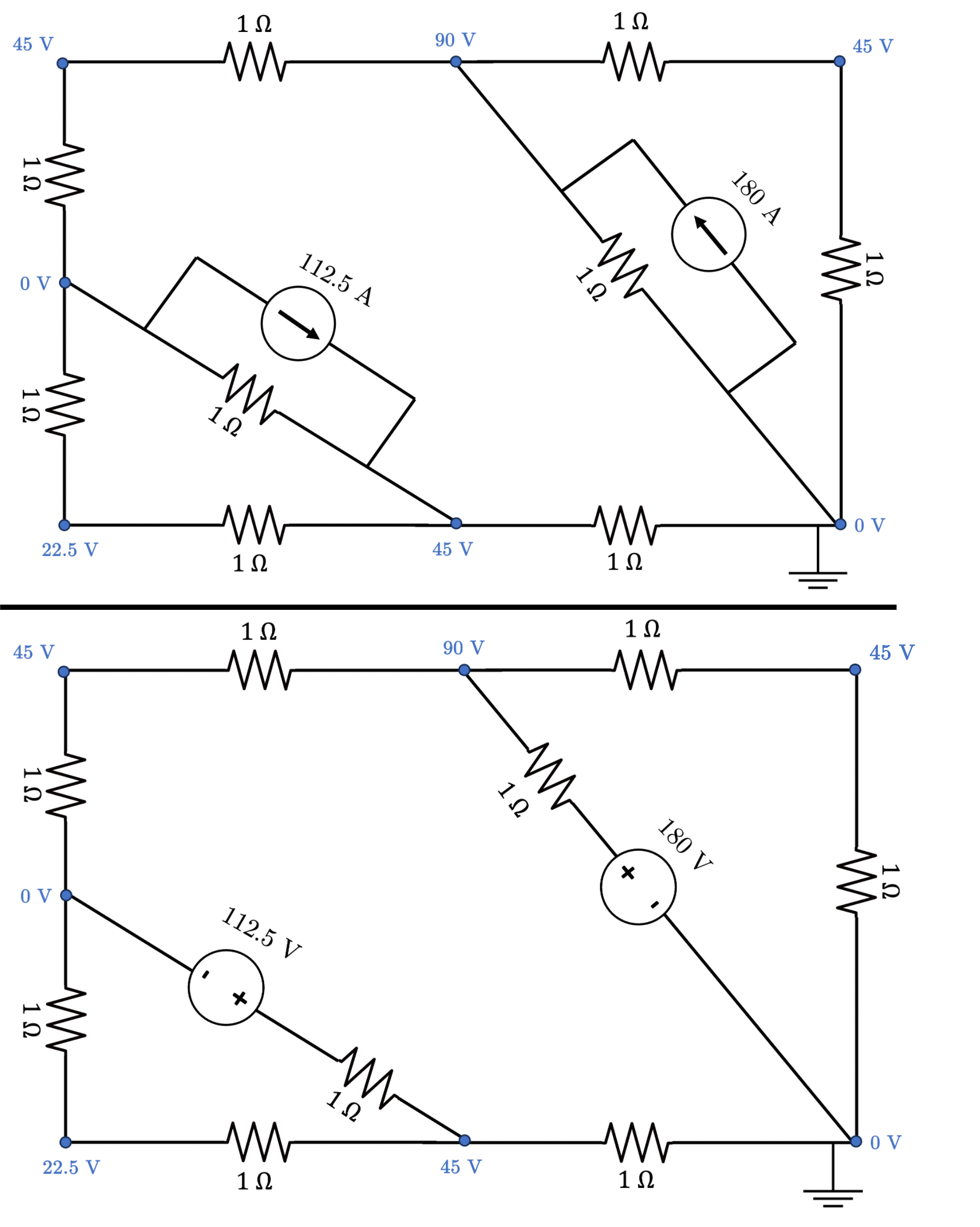}
    \caption{7-bus network from Fig. \ref{fig:7bus} and its DC circuit analogy (top: with current sources representing congestion; bottom: with voltage sources).}
    \label{fig:7bus_total}
\end{figure}

Using the conversions outlined in Sec. \ref{sec:analogy}, the DC circuit equivalent of Fig. \ref{fig:7bus} is shown with both current sources (top) and voltage sources (bottom) in Fig. \ref{fig:7bus_total}. From the ``KCL'' optimality condition \eqref{eqn:dualFeasibility_theta}, it can be seen that the $\bB^T \blambda$ term is analogous to $V/R$ (voltage divided by resistance, or current), and therefore the the $\bD^T \bmu$ term can be thought of as current injections. $\bD$ is unitless, and thus $\bmu$ must be in units of Amps. This is why the current source's value is equal to the original value of the congestion cost $\mu$.

However, if working with voltage sources is preferable, a simple source transformation can be performed on current sources. A source transformation uses Ohm's law to replace a current source in parallel with a resistor with a voltage source in series with a resistor (or vice versa). The resulting value of the voltage source can be found with Ohm's law, where $V = IR$. In this particular example, $R =1 \Omega$, which results in the voltage source having the same magnitude as the current source. The resistor's value does not change in the transformation, but simply moves from a parallel connection to a series connection with the source.

Ground is placed at the location of the cheapest marginal generator, which is at bus 2 in the original network, and the bottom right node in the circuit. Because all of the line susceptances in the original network were 1 pu, all of the resistances in the circuit become 1 $\Omega$. 

\section{Utility of Circuit Interpretations}
After the optimal solution of the DC OPF problem is found and the equivalent circuit is formed, intuitive conclusions can be drawn from the use of DC circuit techniques. Here, we use a few examples of insights that can be gleaned using some of these techniques.

\subsection{Superposition for Congestion Analysis}
The superposition principal, as applied to linear DC circuits, states that any voltage or current throughout the network can be derived by summing the contributions of individual voltage and current sources. In this analogy, voltage and current sources represent congested lines; thus, we can analyze the impact on individual congested lines on prices network-wide. These insights could provide useful for objectives such as transmission planning, for example, where ISOs may desire to upgrade the transmission lines that have the most potential to reduce system-wide prices.

\begin{figure}[t!]
    \centering
    \includegraphics[width=0.8\columnwidth]{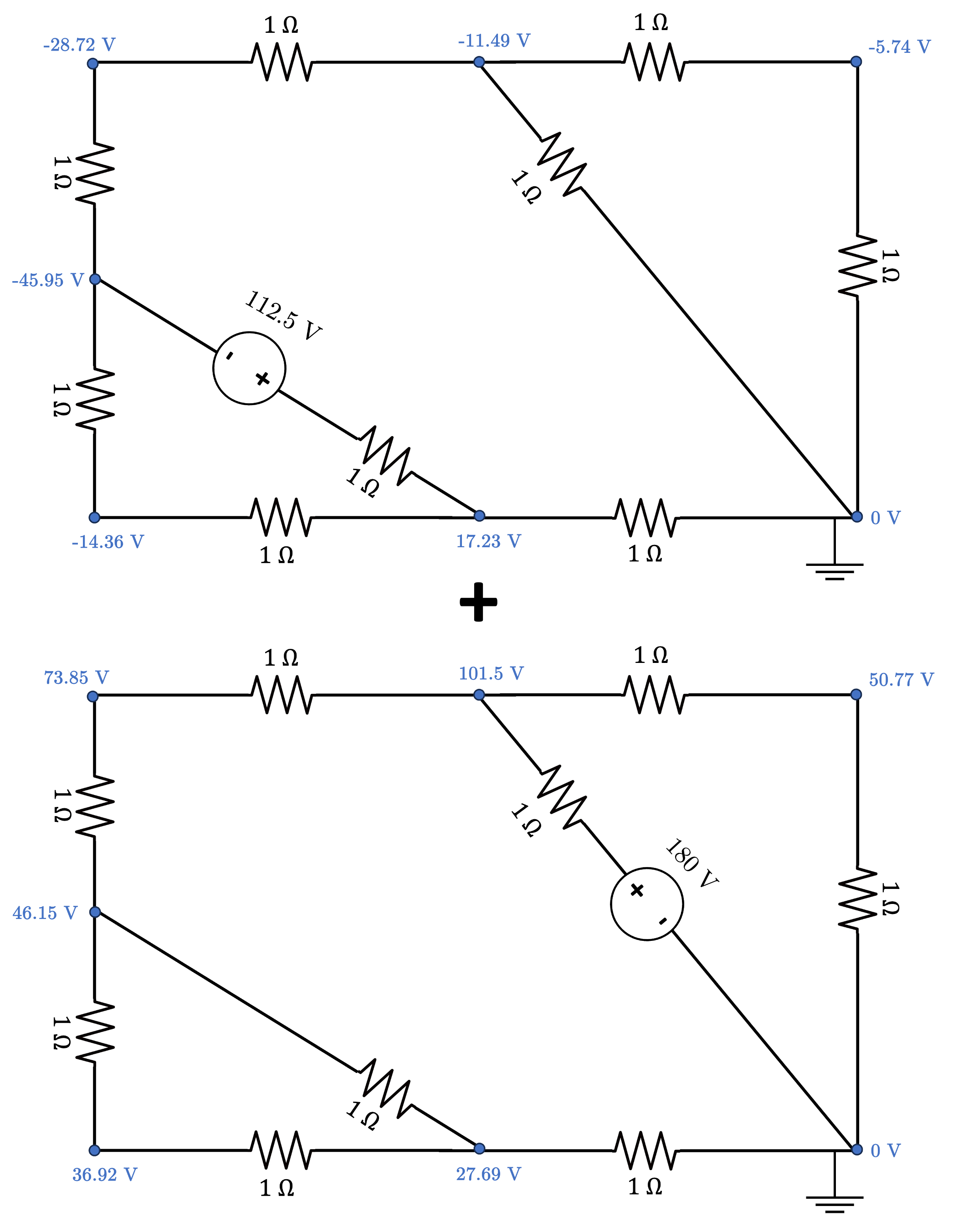}
    \caption{7-bus network with nodal voltage contributions from each source (congested line). The resulting nodal voltages (LMPs) are the sum of the contribution of the individual sources. Note that the sum may be off by a few decimal places when compared to the original voltages shown in Fig. \ref{fig:7bus_total} due to rounding.}
    \label{fig:superposition}
\end{figure}

Using the same 7-bus example as before, we can first assess the impact of the first congested line (112.5 V source), the impact of the second congested line (180 V source), and sum them to recover the original nodal voltages. The impact of each source is shown in Fig. \ref{fig:superposition}, where the source that is not considered is replaced with a shorted wire (current sources would be replaced with open circuits). Note that there may be slight differences to the LMPs in Fig. \ref{fig:7bus_total} due to the rounding / numerical calculations of the circuit solver used. In this example, the contribution of the 112.5 V source is interesting - it actually contributes negative nodal voltages to most of the nodes.

\subsection{Existence of Negative Prices/Voltages}
In general, even without renewable energy generation, negative bids, ``must take'' generation, ramp limits, or minimum on/off times, negative prices can still occur within a network simply due to congestion. Negative nodal voltages occur in a DC circuit when ground is not placed at the node with the lowest potential. Since ground is the reference point that all other nodal voltages are with respect to, if ground (0 V) is located at the lowest potential node, no other node can be lower than it (i.e., negative) by definition.

\begin{figure}[t!]
    \centering
    \includegraphics[width=1\columnwidth]{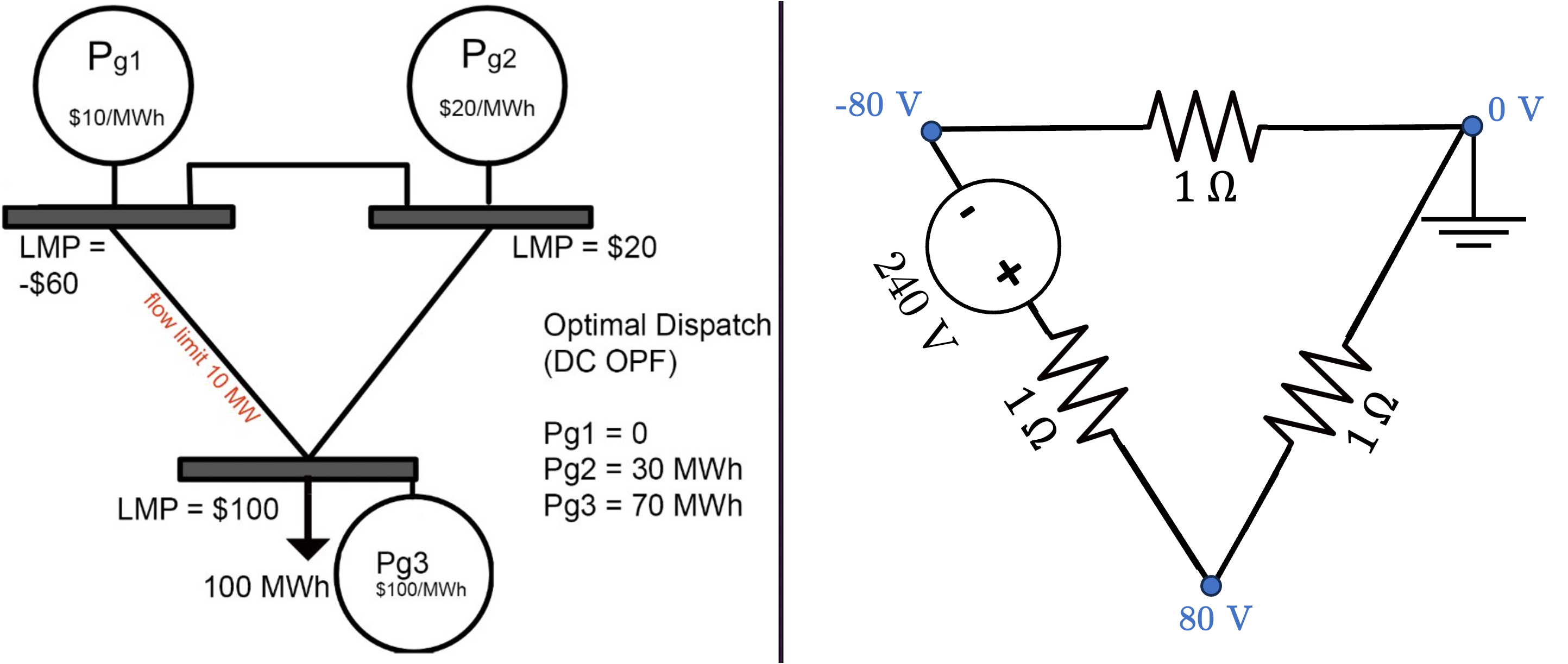}
    \caption{3-bus network with negative prices due to one congested line (left), and corresponding DC circuit equivalent (right). Note that because ground (the cheapest marginal generator) is not located at the node with the lowest potential, negative voltages (LMPs) necessarily result.}
    \label{fig:3bus_neg}
\end{figure}

To illustrate this, consider the simple 3-bus system in Fig. \ref{fig:3bus_neg}. In this network, no upper limits on generation exist, but due to a 10 MW flow constraint between buses 1 and 3, the cheapest generator, generator 1, cannot generate any power to satisfy the load at bus 3 without violating DC power flow constraints. Thus, this generator (1) is not marginal. Generator 2, twice as expensive, becomes the cheapest marginal generator, and ground is placed at bus 2 in the circuit equivalent. Because ground is not located at the node with the lowest potential (here, the negative terminal of the voltage source), negative nodal voltages occur.

Recall that the value of the cheapest marginal generator (\$20/MWh) must be added to each nodal voltage to recover the original LMPs. In this case, 20 is added to the -80 V at bus 1 to recover the original -\$60 LMP at bus 1. In general, we can state the following regarding negative prices: In a circuit, if ground is not placed at the node with the lowest potential, negative voltages occur. In a power network, if the cheapest marginal generator is not located at the bus with the lowest LMP in the network, negative prices are guaranteed to occur. This means that in order to know whether or not there will be negative prices in a network, you simply have to know if the location of the cheapest marginal generator is also the location of the lowest LMP.

\subsection{Uncovering LMPs from Limited Information}
Another application of the DC circuit analogy is that all voltage drops (LMP differences) throughout the network can be recovered by knowing the topology of the network, the resistance (susceptance) of the lines, and the location and value of voltage/current sources (congested lines). This is due to the fact that a DC circuit can be fully solved (e.g. all voltage drops and currents determined) if this information is known. If the location of the cheapest marginal generator is additionally known, one can recover the exact LMPs at each bus; otherwise, only LMP differences can be recovered. This finding can also be derived from the optimality condition \eqref{eqn:dualFeasibility_theta}; Lagrange multipliers corresponding to generator upper/lower limits ($\bgamma$), for example, or any values for primal variables, have to be known to recover LMP differences.

\section{Conclusion}
In this paper, we introduced the concept of modeling locational marginal prices as physical quantities (voltages) within a DC circuit with independent sources. This modeling technique allows for the following benefits: 1) A method to understand and interpret prices without having to understand more advanced mathematical concepts such as Lagrangian duality; 2) Use of widely-known DC circuit analysis techniques such as superposition to analyze the impact of individual congested lines on LMPs; 3) Improved intuition for how LMP differences propagate throughout a network akin to physical ``flows.'' Future work will leverage the concepts introduced in this paper to derive formal mathematical properties of LMPs with respect to network parameters such as topology, distance between marginal generators, and line susceptances.

\bibliographystyle{IEEEtran}
\bibliography{references}

\end{document}